\documentclass[12pt,leqno,draft]{article}
\usepackage{amsfonts}
\usepackage{amsmath, amsthm, amsfonts, amssymb, color}
\usepackage{mathrsfs}
\usepackage{color}
\newtheorem{thm}{Theorem}[section]

\newtheorem{lem}[thm]{Lemma}

\theoremstyle{definition}

\newcommand{\scr}[1]{\mathscr #1}
\definecolor{wco}{rgb}{0.5,0.2,0.3}

\numberwithin{equation}{section} \theoremstyle{remark}

\newcommand{\ua}{\uparrow}

\title{{\bf McKean-Vlasov   SDEs with Drifts  Discontinuous under  Wasserstein Distance}\footnote{Supported in
 part by  NNSFC (11771326, 11831014, 11801406).} }
\author{
{\bf   Xing Huang $^{a)}$,  Feng-Yu Wang $^{a), b)}$  }\\
\footnotesize{ a)Center for Applied Mathematics, Tianjin
University, Tianjin 300072, China}\\
\footnotesize{  xinghuang@tju.edu.cn}\\
 \footnotesize{ b)Department of Mathematics,
Swansea University, Singleton Park, SA2 8PP, United Kingdom}\\
\footnotesize{  wangfy@tju.edu.cn}}
\begin{document}
\allowdisplaybreaks
\def\R{\mathbb R}  \def\ff{\frac} \def\ss{\sqrt} \def\B{\mathbf
B} \def\W{\mathbb W}
\def\N{\mathbb N} \def\kk{\kappa} \def\m{{\bf m}}
\def\ee{\varepsilon}\def\ddd{D^*}
\def\dd{\delta} \def\DD{\Delta} \def\vv{\varepsilon} \def\rr{\rho}
\def\<{\langle} \def\>{\rangle} \def\GG{\Gamma} \def\gg{\gamma}
  \def\nn{\nabla} \def\pp{\partial} \def\E{\mathbb E}
\def\d{\text{\rm{d}}} \def\bb{\beta} \def\aa{\alpha} \def\D{\scr D}
  \def\si{\sigma} \def\ess{\text{\rm{ess}}}
\def\beg{\begin} \def\beq{\begin{equation}}  \def\F{\scr F}
\def\Ric{\text{\rm{Ric}}} \def\Hess{\text{\rm{Hess}}}
\def\e{\text{\rm{e}}} \def\ua{\underline a} \def\OO{\Omega}  \def\oo{\omega}
 \def\tt{\tilde} \def\Ric{\text{\rm{Ric}}}
\def\cut{\text{\rm{cut}}} \def\P{\mathbb P} \def\ifn{I_n(f^{\bigotimes n})}
\def\C{\scr C}      \def\aaa{\mathbf{r}}     \def\r{r}
\def\gap{\text{\rm{gap}}} \def\prr{\pi_{{\bf m},\varrho}}  \def\r{\mathbf r}
\def\Z{\mathbb Z} \def\vrr{\varrho} \def\ll{\lambda}
\def\L{\scr L}\def\Tt{\tt} \def\TT{\tt}\def\II{\mathbb I}
\def\i{{\rm in}}\def\Sect{{\rm Sect}}  \def\H{\mathbb H}
\def\M{\scr M}\def\Q{\mathbb Q} \def\texto{\text{o}} \def\LL{\Lambda}
\def\Rank{{\rm Rank}} \def\B{\scr B} \def\i{{\rm i}} \def\HR{\hat{\R}^d}
\def\to{\rightarrow}\def\l{\ell}\def\iint{\int}
\def\EE{\scr E}\def\Cut{{\rm Cut}}
\def\A{\scr A} \def\Lip{{\rm Lip}}
\def\BB{\scr B}\def\Ent{{\rm Ent}}\def\L{\scr L}
\def\R{\mathbb R}  \def\ff{\frac} \def\ss{\sqrt} \def\B{\mathbf
B}
\def\N{\mathbb N} \def\kk{\kappa} \def\m{{\bf m}}
\def\dd{\delta} \def\DD{\Delta} \def\vv{\varepsilon} \def\rr{\rho}
\def\<{\langle} \def\>{\rangle} \def\GG{\Gamma} \def\gg{\gamma}
  \def\nn{\nabla} \def\pp{\partial} \def\E{\mathbb E}
\def\d{\text{\rm{d}}} \def\bb{\beta} \def\aa{\alpha} \def\D{\scr D}
  \def\si{\sigma} \def\ess{\text{\rm{ess}}}
\def\beg{\begin} \def\beq{\begin{equation}}  \def\F{\scr F}
\def\Ric{\text{\rm{Ric}}} \def\Hess{\text{\rm{Hess}}}
\def\e{\text{\rm{e}}} \def\ua{\underline a} \def\OO{\Omega}  \def\oo{\omega}
 \def\tt{\tilde} \def\Ric{\text{\rm{Ric}}}
\def\cut{\text{\rm{cut}}} \def\P{\mathbb P} \def\ifn{I_n(f^{\bigotimes n})}
\def\C{\scr C}      \def\aaa{\mathbf{r}}     \def\r{r}
\def\gap{\text{\rm{gap}}} \def\prr{\pi_{{\bf m},\varrho}}  \def\r{\mathbf r}
\def\Z{\mathbb Z} \def\vrr{\varrho} \def\ll{\lambda}
\def\L{\scr L}\def\Tt{\tt} \def\TT{\tt}\def\II{\mathbb I}
\def\i{{\rm in}}\def\Sect{{\rm Sect}}  \def\H{\mathbb H}
\def\M{\scr M}\def\Q{\mathbb Q} \def\texto{\text{o}} \def\LL{\Lambda}
\def\Rank{{\rm Rank}} \def\B{\scr B} \def\i{{\rm i}} \def\HR{\hat{\R}^d}
\def\to{\rightarrow}\def\l{\ell}
\def\8{\infty}\def\I{1}\def\U{\scr U}
\maketitle

\begin{abstract} Existence and uniqueness are proved for Mckean-Vlasov type distribution dependent SDEs with singular drifts satisfying an integrability condition in space variable and the Lipschitz condition in distribution variable with respect to $\W_0$ or $\W_0+\W_\theta$ for some $\theta\ge 1$, where $\W_0$ is the total variation distance and $\W_\theta$ is the $L^\theta$-Wasserstein distance.
This improves some existing results where the drift is continuous in  the distribution variable  with respect to the Wasserstein distance.

\end{abstract} \noindent
 AMS subject Classification:\  60H1075, 60G44.   \\
\noindent
 Keywords: Distribution dependent SDEs, total variation distance, Wasserstein distance, Krylov's estimate, Zvonkin's transform.
 \vskip 2cm

\section{Introduction}

 Consider the following   distribution dependent SDE on $\R^d$:
\beq\label{E1} \d X_t= b_t(X_t, \L_{X_t})\d t +\si_t(X_t,\L_{X_t})\d W_t,\ \ \ t\in [0,T],\end{equation}
where $T>0$ is a fixed time, $(W_t)_{t\in [0,T]}$ is the $m$-dimensional Brownian motion on a complete filtration probability space $(\OO,\{\F_t\}_{t\in[0,T]},\P)$, $\L_{X_t}$ is the law of $X_t$,
$$b: [0,T]\times\R^d\times \scr P\to \R^d,\ \ \si: [0,T]\times\R^d\times\scr P\to \R^d\otimes\R^m$$ are measurable, and $\scr P$ is the space of all probability measures on $\R^d$ equipped with the weak topology.

This type SDEs are also called   McKean-Vlasov   SDEs   and  mean field SDEs, and have been intensively investigated due to its wide applications, see for instance
  \cite{BR1,BR2, CA, CF, DV1,DV2, HRW,  SZ, FYW1} and references within.

An adapted continuous process on $\R^d$ is called a (strong) solution of \eqref{E1}, if
\beq\label{A1} \E\int_0^T\big\{|b_t(X_t,\L_{X_t})|+\|\si_t(X_t,\L_{X_t})\|^2\big\}\d t<\infty,\end{equation}
and $\P$-a.s.
\beq\label{A2} X_t= X_0 +\int_0^t b_s(X_s,\L_{X_s})\d s +\int_0^t \si_s(X_s,\L_{X_s})\d W_s,\ \ t\in [0,T].\end{equation} We call \eqref{E1} (strongly) well-posed for an $\F_0$-measurable  initial value $X_0$, if \eqref{E1} has a unique solution starting at $X_0$.

When a different probability measure $\tt\P$ is concerned, we use $\L_\xi|\tt \P$ to denote the law of a random variable $\xi$ under the probability $\tt\P$, and use $\E_{\tt\P}$ to stand for  the expectation under $\tt\P$.
For any $\mu_0\in \scr P$,   $(\tt X_t,\tt W_t)_{t\in [0,T]}$  is called a weak solution to \eqref{E1} starting at $\mu_0$, if
$(\tt W_t)_{t\in [0,T]}$ is the $m$-dimensional Brownian motion under a   complete filtration probability space $(\tt\OO,\{\tt\F_t\}_{t\in[0,T]},\tt\P),$ $(\tilde{X}_t)_{t\in [0,T]}$ is   a continuous $\tt\F_t$-adapted process on $\R^d$ with $\L_{\tt X_0}|\tt \P=\mu_0$, and
\eqref{A1}-\eqref{A2} hold for $(\tt X, \tt W, \tt\P,\E_{\tt P})$ replacing  $(X, W, \P,\E).$
 We call \eqref{E1} weakly well-posed for an initial distribution $\mu_0$, if it has a unique weak solution starting at $\mu_0$; i.e. it has a weak solution $(\tt X_t,\tt W_t)_{t\in [0,T]}$  with initial distribution $\mu_0$ under some complete filtration probability space    $(\tt\OO,\{\tt\F_t\}_{t\in[0,T]},\tt\P)$, and $\L_{\tt X_{[0,T]}}|\tt \P= \L_{\bar X_{[0,T]}}|\bar \P$ holds for any other weak solution with the same initial distribution
$(\bar X_t,\bar W_t)_{t\in [0,T]}$  under some complete filtration probability space $(\bar\OO,\{\bar\F_t\}_{t\in[0,T]},\bar\P).$

Recently, the (weak and strong) well-posedness  is studied  in  \cite{BB,BBP,CR,HW,L,MV,RZ}  for \eqref{E1}   with $\si_t(x,\gamma)=\si_t(x)$ independent of the distribution variable $\gamma$, and with singular drift $b_t(x,\gamma)$. See also \cite{HRW,L} for the case with memory.
We briefly recall some  conditions on $b$ which together with a regular and non-degenerate condition on $\si$ implies the well-posedness of \eqref{E1}.  To this end, we recall the $L^\theta$-Wasserstein distance $\W_\theta$
for   $\theta>0$:
  $$\W_\theta(\gamma,\tilde{\gamma}):= \inf_{\pi\in \C(\gamma,\tilde{\gamma})} \bigg(\int_{\R^d\times\R^d} |x-y|^\theta \pi(\d x,\d y)\bigg)^{\ff 1 {1\lor\theta}},\ \ \gamma,\tilde{\gamma}\in \scr P,$$
  where $\C(\gamma,\tilde{\gamma})$ is the set of all couplings of $\gamma$ and $\tilde{\gamma}$.
By the convention that $r^0= 1_{\{r>0\}}$ for $r\ge 0$, we may regard $\W_0$ as the total variation distance, i.e. set
$$\W_0(\gamma,\tilde{\gamma})= \|\gamma-\tilde{\gamma}\|_{TV}:=\sup_{A\in \B(\R^d)} |\gamma(A)-\tilde{\gamma}(A)|.$$

 References   \cite{BB,BBP} give  the well-posedness of  \eqref{E1} with a deterministic initial value $X_0\in \R^d$, where the drift $b_t(x,\gamma)$ is assumed to be linear growth in $x$ uniformly in $t,\gamma$,
and
$$|b_t(x,\gamma)-b_t(x,\tilde{\gamma})|\le \phi(\W_1(\gamma,\tilde{\gamma}))$$ holds  for some function $\phi\in C((0,\infty);(0,\infty))$  with $\int_0^\cdot\ff 1 {\phi(s)}\d s=\infty$.
Note that for distribution dependent SDEs the well-posedness for deterministic initial values does not imply that for random ones.

  \cite[Theorem 3]{MV}  presents the well-posedness of \eqref{E1} with exponentially integrable $X_0$ and
a drift  $b$   of type
\beq\label{B*} b_t(x,\gamma):= \int_{\R^d} \tt b_t(x,y)\gamma(\d y),\end{equation}
where $\tt b_t(x,y)$ has linear growth in $x$ uniformly in $t$ and $y$. Since $\tt b_t(x,y)$ is bounded in $y$, $b_t(x,\cdot)$ is Lipschtiz continuous in the total variation distance $\W_0$.
\cite{RZ} considers the same type drift and proves the well-posedness  of \eqref{E1}  under the conditions that  $\E|X_0|^\bb<\infty$ for some $\bb>0$ and
$$|\tt b_t(x,y)|\le h_t(x-y)$$ for some $h\in L^q([0,T];\tt L^p(\R^d))$ for   some $p,q>1$ with $\ff d p+\ff 2 q<1$, where $\tt L^p$ is a localized $L^p$ space.

In \cite{CR}  the well-posedness of \eqref{E1} is proved for $X_0$ satisfying $\E|X_0|^2<\infty$, and  for $b$   given by
\beq\label{*RP} b_t(x,\gamma)= \tt b_t(x, \gamma(\varphi)),\end{equation}
where $\gamma(\varphi):=\int_{\R^d} \varphi\d\gamma$ for some $\aa$-H\"older  continuous function $\varphi$,  and $|\tt b_t(x,r)|+ | \pp_r \tt b_t(x,r)|$ is bounded.  Consequently, $b_t(x,\gamma)$ is bounded and Lipschitz continuous in $\gamma$ with respect to $\W_\aa$.

In \cite{HW} the well-posedness is derived under the conditions that    $\E|X_0|^\theta<\infty$ for some $\theta\ge 1$,   $b_t(x,\gamma)$ is Lipschitz continuous   in $\gamma$ with respect to $\W_\theta$,
and   for any $\mu\in C([0,T];\scr P_\theta),$
$$b^\mu_t(x):= b_t(x,\mu_t),\ \ (t,x)\in [0,T]\times \R^d$$ satisfies  $|b^\mu|^2\in L_{p}^{q}(T)$ for some $(p,q) \in \scr K,$ where
\beg{align*}  &L_p^q(T):=\bigg\{f\in\B([0,T]\times\R^d):  \int_0^T \Big(\int_{\R^d} |f_t(x)|^p\d x\Big)^{\ff q p} \d t <\infty\bigg\},\\
&\scr K:=\Big\{(p,q)\in (1,\infty)\times(1,\infty):\   \ff d p +\ff 2 q<2\Big\}.\end{align*}

Moreover, in \cite{JM} the well-posedness of \eqref{E1} has been proved for
\beq\label{*PP} b_t(x,\mu)=\tt b(\rr_\mu(x)),\ \ \si_t(x,\mu)=\tt \si(\rr_\mu(x))\end{equation}
with initial distribution having density function (with respect to the Lebesgue mseaure) in the class $H^{2+
\aa}$ for some $\aa>0$,
where $\rr_\mu$ is the density function of $\mu$ with respect to the Lebesgue measure, $\tt b\in C^2([0,\infty);\R^d)$ and $\tt\si\in C^3([0,\infty);\R^d\otimes\R^d).$  As for the weak well-posedness,   \cite{J0} assumes that  $b$ is bounded and $\W_0$-Lipschitz continuous in distribution variable, and $\sigma$ is Lipschitz continuous in space
variable.

\

 In this paper,   we  prove the (weak and strong) well-posedness  of  \eqref{E1}  for   general type $b$ with $b_t(x,\gamma)$  Lipschitz continuous in $\gamma$ under the metric $\W_0$ or $\W_0+\W_\theta$ for some $\theta\ge 1$.
 This condition is weaker than those in \cite{BB,BBP,CR,HW} in the sense that the drift is not necessarily continuous in the Wasserstein distance, but is incomparable with those in \cite{MV,RZ} where  $b$ is of the integral type as in \eqref{B*}. Moreover, our result works for any initial value and initial distribution.

Recall that a continuous function $f$ on $\mathbb{R}^d$ is called weakly differentiable, if there exists (hence unique) $\xi\in L^1_{loc}(\mathbb{R}^d)$ such that
$$\int_{\mathbb{R}^d}(f\Delta g)(x)\d x=-\int_{\mathbb{R}^d}\langle \xi,\nabla g\rangle (x)\d x, \ \ g\in C_0^\infty(\mathbb{R}^d).$$
In this case, we write $\xi=\nabla f$ and call it the weak gradient of $f$. For $p,q\ge 1$,  let
$$L_{p,loc}^q(T)=\bigg\{f\in \B([0,T]\times\R^d): \int_0^T  \Big(\int_K |f_t(x)|^p\d x\Big)^{\ff q p} \d t<\infty,\ K\subset \R^d \ \text{\ is \ compact}\bigg\}.$$
We will use the following conditions.
\begin{enumerate}
\item[$(A_\si)$] $\si_t(x,\gamma)=\si_t(x)$ is uniformly continuous in $x\in\R^d$ uniformly in $t\in [0,T];$ the weak gradient $\nn \si_t$ exists for a.e. $t\in [0,T]$ such that
    $|\nn \si|^2\in L_p^q(T)$ for some  $(p,q)\in \scr K$;  and there exists a constant $K_1\ge 1$ such that
\begin{align} \label{si-} K_1^{-1} I \le (\sigma_t\sigma_t^\ast)(x)\le K_1 I,\ \ (t,x)\in [0,T]\times\R^d,\end{align} where $I$ is the $d\times d$ identity matrix.
 \item[$(A_b)$]
$b=\bar{b}+\hat{b}$, where $\bar b$ and $\hat b$ satisfy
\begin{equation}\label{con}\beg{split}
&|\hat{b}_t(x,\gamma)-\hat{b}_t(y,\tilde{\gamma})|+|\bar{b}_t(x,\gamma)-\bar{b}_t(x,\tilde{\gamma})|\\
&\le K_2(\|\gamma-\tilde{\gamma}\|_{TV}+\W_\theta(\gamma,\tilde{\gamma})+|x-y|),\ \ t\in[0,T], x,y \in \R^d, \gamma,\tilde{\gamma}\in \scr P_\theta\end{split}
\end{equation}for some constants $\theta,K_2\ge 1,$ and
 there exists $(p,q)\in \scr K$ such that
\beq\label{ihp} \sup_{t\in [0,T]} |\hat{b}_t(0,\delta_0)|+\sup_{\mu\in C([0,T]; \scr P_{\theta})}\| |\bar{b}^\mu|^2\|_{ L_{p}^{q}(T)}<\infty,\end{equation}
where $\bar b_t^\mu(x):= \bar b_t(x,\mu_t)$ for $(t,x)\in [0,T]\times\R^d$, and  $\delta_0$ stands for the Dirac measure at the  point $0\in\R^d$.
\item[$(A_b')$] For any $\mu\in \B([0,T];\scr P)$, $|b^\mu|^2\in L_{p,loc}^q(T)$ for some $(p,q)\in \scr K$. Moreover, there exists a function $\Gamma:[0,\infty)\to [0,\infty)$ satisfying $\int_{1}^\infty\frac{1}{\Gamma(x)}=\infty$ such that
\begin{align}\label{ne}\<b_t(x,\delta_0),x\>\leq \Gamma(|x|^2),\ \ t\in[0,T],x\in \R^d.\end{align}
In addition, there exists a constant $K_3\ge 1$ such that
 \begin{equation}\label{cono}\beg{split}
&|b_t(x,\gamma)-b_t(x,\tilde{\gamma})|\le K_3\|\gamma-\tilde{\gamma}\|_{TV},\ \ t\in[0,T], x\in \R^d, \gamma,\tilde{\gamma}\in \scr P.
\end{split}
\end{equation}
\end{enumerate}
When \eqref{E1} is weakly well-posed for initial distribution $\gamma$, we denote $P_t^*\gamma$ the distribution of the weak solution at time $t$.
\begin{thm}\label{T1.1} Assume $(A_\si)$.
\beg{enumerate}
\item[$(1)$] If $(A_b')$ holds, then  \eqref{E1}  is strongly and weakly well-posed for any initial values and any initial distribution. Moreover,
\beq\label{B1} \|P_t^*\mu_0- P_t^*\nu_0\|_{TV}^2\le 2\e^{\ff{K_1K_3^2t}2}\|\mu_0-\nu_0\|_{TV}^2,\ \ t\in [0,T], \mu_0,\nu_0\in \scr P.\end{equation}
 \item[$(2)$]   Let $\E|X_0|^\theta<\infty$ and $\mu_0(|\cdot|^\theta)<\infty$.  If  $(A_b)$ holds,   then  $\eqref{E1}$ is strongly well-posed  for initial value $X_0$ and weakly well-posed for initial distribution $\mu_0$.  Moreover,
 there exists a constant $c>0$ such that for any $\mu_0,\nu_0\in \scr P_\theta,$
\beq\label{B2}\begin{split} &\|P_t^*\mu_0- P_t^*\nu_0\|_{TV} +\W_\theta(P_t^*\mu_0, P_t^*\nu_0)\\
&\le c \big\{\|\mu_0-\nu_0\|_{TV}+W_\theta(\mu_0,\nu_0)\big\}, \ \ t\in [0,T].\end{split}\end{equation}
 \end{enumerate}
\end{thm}

To illustrate this result comparing with earlier ones, we present an example of $b$  which satisfies our conditions but is not of  type \eqref{B*}-\eqref{*PP} and is discontinuous in both the space variable 
 and the distribution variable under the weak topology. If one wants to control a stochastic system in terms of an ideal reference distribution $\mu_0$, it is natural to take a drift depending on a probability distance between $\mu_0$ and
the law of the system. As two typical probability distances, the total variation and Wasserstein distances have been widely applied in applications. So, we take for instance
$$b_t(x,\mu)= \bar b(t,x,\mu)+ h(t, x, \W_\theta(\mu,\mu_0), \|\mu-\mu_0\|_{TV})$$
for some $\theta\ge 1,$ where $\bar b$ satisfies \eqref{con} and \eqref{ihp} for $\hat b=0$ which refers to the singularity in the space variable $x$,
and $h: [0,T]\times \R^d\times [0,\infty)^2\to\R^d$ is measurable such that $h(t,x,r,s)$ is bounded in $t\in [0,T]$ and Lipschitz continuous in $(x,r,s)\in\R^d\times [0,\infty)^2$ uniformly in $t\in [0,T].$   Obviously, $b(t,x,\mu)$  satisfies condition $(A_b)$ but is not of type \eqref{B*}-\eqref{*PP} and can be discontinuous in $x$ and  $\mu$ under the weak topology. 

In the next section we make some preparations, which will be used in Section 3 for the proof of Theorem \ref{T1.1}.

  \section{Preparations}

We first present the following version of  Yamada-Watanabe principle modified from   \cite[Lemma 3.4]{HW}.

\beg{lem}\label{YW}  Assume that  $\eqref{E1}$ has a weak solution $(\bar X_t)_{t\in [0,T]}$ under probability $\bar \P$, and let $\mu_t=\L_{\bar X_t}|\bar \P, t\in [0,T].$
  If the SDE
\beq\label{CSDE} \d X_t= b_t(X_t,\mu_t)\d t+ \si_t(X_t,\mu_t)\d W_t\end{equation}has strong uniqueness for some initial value $X_0$ with $\L_{X_0}=\mu_0$, then $\eqref{E1} $ has a strong solution starting at $X_0$.
If moreover \eqref{E1} has strong uniqueness for any initial value $X_0$ with $\L_{X_0}=\mu_0$, then it is  weakly well-posed  for the initial distribution $\mu_0$.
\end{lem}
\beg{proof}(a) Strong existence. Since $\mu_t=\L_{\bar X_t}|\bar \P$, $\bar X_t$ under $\bar \P$ is also a weak solution of \eqref{CSDE} with initial distribution $\mu_0$.
By the  Yamada-Watanabe principle, the strong uniqueness of \eqref{CSDE} with initial value $X_0$ implies the strong (resp. weak) well-posedness of \eqref{CSDE} starting at $X_0$ (resp. $\mu_0)$. In particular, the weak uniqueness implies $\L_{X_t}=\mu_t, t\in [0,T]$, so that $X_t$ solves  \eqref{E1}.

(b) Weak uniqueness. Let $\tt X_t$ under probability $\tt \P$ be another weak solution of $\eqref{E1}$ with initial distribution $\mu_0$. For any initial value $X_0$ with $\L_{X_0}=\mu_0$,
the strong uniqueness of \eqref{CSDE} starting at $X_0$ implies
$$X_{[0,T]}=F(X_0,W_{[0,T]})$$ for some measurable function $F: \R^d\times C([0,T];\R^d)\to C([0,T];\R^d).$ This and the weak uniqueness of \eqref{CSDE} proved in (a) yield
\beq\label{PPO} \L_{\bar X_{[0,T]}}|\bar\P= \L_{ X_{[0,T]}}|\P.\end{equation}
Let $\hat X_{[0,T]}=F(\tt X_0,\tt W_{[0,T]})$. We have  $\hat X_0= \tt X_0$ and
$$\L_{\hat X_{[0,T]}}|\tt \P= \L_{X_{[0,T]}}|\P.$$
This and \eqref{PPO} imply  $\L_{\hat X_t}|\tt\P=\mu_t$, so that $\hat X_t$ under $\tt\P$ is a weak solution of \eqref{E1} with $\hat X_0=\tt X_0$.
By the strong uniqueness of \eqref{E1}, we derive $\hat X_{[0,T]}=\tt X_{[0,T]}$. Combining this with \eqref{PPO} we obtain
$$\L_{\tt X_{[0,T]}}|\tt\P= \L_{\hat X_{[0,T]}}|\tt\P= \L_{ X_{[0,T]}}|\P= \L_{\bar X_{[0,T]}}|\bar\P,$$ i.e.    \eqref{E1} has weak uniqueness starting at $\mu_0$.
\end{proof}
We will use the following result for the maximal operator:
\begin{align}\label{max}
\M h(x):=\sup_{r>0}\frac{1}{|B(x,r)|}\int_{B(x,r)}h(y)\d y,\ \  h\in L^1_{loc}(\mathbb{R}^d), x\in \R^d,
 \end{align} where $B(x,r):=\{y: |x-y|<r\},$  see   \cite[Appendix A]{CD}.

 \begin{lem} \label{Hardy}  There exists a constant $C>0$ such that for any continuous and weak differentiable function $f$,
 \beq\label{HH1}
|f(x)-f(y)|\leq C|x-y|(\M |\nabla f|(x)+\M |\nabla f|(y)),\ \  {\rm a.e.}\ x,y\in\R^d.\end{equation}
Moreover, for any $p>1$,    there exists a constant $C_{p}>0$ such that
\beq\label{HH2}
\|\M f\|_{L^p}\leq C_{p}\|f\|_{L^p},\ \ f\in L^p(\R^d).
 \end{equation}
\end{lem}

To compare the distribution dependent SDE \eqref{E1} with a classical one, for any $\mu\in \B([0,T]; \scr P),$ let $b_t^\mu(x):=b_t(x,\mu_t)$ and consider the classical SDE
\beq\label{EW1} \d X^\mu_t= b_t^\mu(X^\mu_t)\d t +\si_t(X^\mu_t)\d W_t,\ \ t\in [0,T].\end{equation}
According to \cite{Z2},  assumption  $(A_\si)$ together with $(A_b)$ or $(A_b')$ implies the strong well-posedness,
where under $(A_b')$ the non-explosion is implied by \eqref{ne}.  For any $\gg\in \scr P$,
 Let $\Phi_t^\gg(\mu)=\L_{X_t^\mu}$ for $(X_t^\mu)_{t\in [0,T]}$ solving \eqref{EW1}  with $\L_{X_0^\mu}=\gg.$
We have the following result.

\beg{lem}\label{LN} Assume $(A_\si)$ and let $\gg\in \scr P$.
\beg{enumerate} \item[$(1)$] If $(A_b')$ holds, then   for any $\mu,\nu\in \B([0,T];\scr P),$
\beq\label{ERT1} \|\Phi_t^\gg(\mu)-\Phi_t^\gg(\nu)\|_{TV}^2\le \ff{K_1K_3^2}4  \int_0^t  \|\mu_s-\nu_s\|_{TV}^2 \d s,\ \ t\in [0,T].\end{equation}
 \item[$(2)$]   If $(A_b)$ holds and $\gg\in \scr P_\theta$, then for any $\mu\in C([0,T]; \scr P_\theta)$, we have $\Phi_\cdot^\gg(\mu)\in C([0,T];\scr P_\theta)$. Moreover, for any $m\ge 1\lor\ff \theta 2,$
  there exists a constant $C>0$ such that for any $\mu,\nu\in C([0,T];\scr P_\theta)$ and $\gg_1,\gg_2\in \scr P_\theta$,
\beq\label{ERT2} \beg{split} & \{ \W_\theta(\Phi_t^{\gg_1}(\mu), \Phi_t^{\gg_2}(\nu))\}^{2m} \\
&\le C\W_\theta(\gg_1,\gg_2)^{2m}+ C \int_0^t \big\{\|\mu_s-\nu_s\|_{TV}+\W_\theta(\mu_s,\nu_s)\big\}^{2m}\d s,\ \ t\in [0,T].\end{split}\end{equation}
\end{enumerate} \end{lem}
 \beg{proof}
(1) Let $(A_b')$ hold and take $\mu,\nu\in \B([0,T];\scr P).$ To compare $\Phi_t^{\gg}(\mu)$ with $\Phi_t^{\gg}(\nu)$, we rewrite \eqref{EW1} as
\beq\label{EW} \d X^\mu_t=b_t (X^\mu_t, \nu_t)\d t+\si_t(X^\mu_t)\d \tt W_t,  \end{equation}
where
$$\tt W_t=W_t+\int_{0}^{t}\xi_s\d s,\ \ \xi_s:=   \{\si_s^\ast(\si_s\si_s^\ast)^{-1}\}(X^\mu_s)[b_s(X^\mu_s, \mu_s)-b_s(X^\mu_s, \nu_s)],\ \ s, t\in [0,T].$$
 Noting that \eqref{si-} together with 
\eqref{cono} implies
 \beq\label{GCC} \E[\e^{\ff 1 2 \int_0^T |\xi_s|^2\d s}]<\infty,\end{equation}
  by the Girsanov theorem we see that $R_T:=\e^{ - \int_0^T \<\xi_s,\d W_s\>-\ff 1 2\int_0^T |\xi_s|^2\d s}$ is a probability density with respect to $\P$, and $(\tt W_t)_{t\in [0,T]}$ is a $d$-dimensional Brownian motion under the probability $\Q:= R_T\P.$

By the weak uniqueness of \eqref{EW1} and $\L_{X_0^\mu}|\Q= \L_{X_0^\mu} =\gg$, we conclude  from \eqref{EW} with  $\Q$-Brownian motion  $\tt W_t$ that
$$\Phi_t^\gg(\nu)= \L_{X_t^\mu }|\Q,\ \ t\in [0,T].$$
Combining this with $(A_\si)$ and applying Pinker's inequality   \cite{Pin}, we obtain
 \beq\label{ESTM}\beg{split}
&2\|\Phi_t^\gg(\nu)-\Phi_t^\gg(\mu)\|_{TV}^2\leq 2 \sup_{\|f\|_\infty\le 1} (\E |f(X_t^\mu) (R_t-1)|)^2
= 2 (\E|R_t-1|)^2 \\
&\leq \E [R_t\log R_t] = \ff 1 2 \E_{\Q}\int_0^t |\xi_s|^2\d s \\
 &\le \ff{K_1} 2 \E_{\Q}\int_0^t |b_s(X^\mu_s, \mu_s)-b_s(X^\mu_s, \nu_s)|^2\d s.\end{split}
 \end{equation}
By  $(A_b')$, this implies   \eqref{ERT1}.

 (2) Let $(A_b)$ hold and take $m\ge 1\lor\ff \theta 2.$ Take $\F_0$-measurable random variables $X_0^\mu$ and $X_0^\nu$ such that $\L_{X_0^\mu}=\gg_1, \L_{X_0^\nu}=\gg_2$ and
 $$\E|X_0^\mu-X_0^\nu|^\theta =\{\W_\theta(\gg_1,\gg_2)\}^\theta.$$
 Let $X_t^\mu$ solve \eqref{EW1} and $X_t^\nu$ solve the same SDE for $\nu$ replacing $\mu$.   We need to find a constant $C>0$ such that for any $t\in [0,T]$,
 \beq\label{ENNB}  \beg{split} &\{\W_\theta(\Phi_t^{\gg_1}(\mu),\Phi_t^{\gg_2}(\nu))\}^{2m}\\
 &\le C(\E|X_0^\mu-X_0^\nu|^\theta)^{\ff {2m}\theta}+  C \int_0^t (\W_\theta(\mu_s,\nu_s)+\|\mu_s-\nu_s\|_{TV})^{2m}\d s,\ \ t\in [0,T].\end{split} \end{equation}
 To this end, we make  a Zvokin type transform as in  \cite{HW} and \cite{SQZ}.

  For any $\lambda>0$, consider the following PDE for $u: [0,T]\times\R^d\to \R^d$:
\beq\label{PDE}
\frac{\partial u_t}{\partial t}+\frac{1}{2}\mathrm{Tr} (\sigma_t\sigma_t^\ast\nabla^2u_t)+\nabla_{b_t^\mu}u_t+\bar{b}_t^\mu  =\lambda u_t,\ \ u_T=0.
\end{equation}
According to \cite[Remark 2.1, Proposition 2.3 (2)]{SQZ}, under assumptions $(A_\si)$ and $(A_b)$,
when $\ll$ is large enough  \eqref{PDE} has a unique solution $\mathbf{u}^{\lambda,\mu}$ satisfying
\begin{align}\label{u0}
\|\mathbf{u}^{\lambda,\mu}\|_\infty+ \|\nabla \mathbf{u}^{\lambda,\mu}\|_{\infty}\leq \frac{1}{5},
\end{align}
and \beq\label{u01} \|\nabla^2 \mathbf{u}^{\lambda,\mu}\|_{L^{2q}_{2p}(T)}<\infty.\end{equation}
Let $\Theta^{\lambda,\mu}_t(x)=x+\mathbf{u}^{\lambda,\mu}_t(x)$. It is easy to see that \eqref{PDE}   and the It\^o formula imply
\beq\label{E-X}
\d \Theta^{\lambda,\mu}_t(X_t^\mu)= (\lambda \mathbf{u}^{\lambda,\mu}_t+\hat{b}_t^\mu) (X_t^\mu)\d t+ (\{\nabla\Theta_t^{\lambda,\mu}\}\sigma_t)(X_t^\mu)\,\d W_t.
\end{equation}
In particular, \eqref{u0} and $\E[|X_0^\mu|^\theta]<\infty$ imply that $\E[|\Theta_0^{\ll,\mu}(X_0^\mu)|^\theta]<\infty$ and  \eqref{E-X}
is an SDE for $\xi_t:= \Theta_t^{\ll,\mu}(X_t^\mu)$ with coefficients of at most linear growth, so that $\L_{\xi_\cdot}\in C([0,T];\scr P_\theta)$ and so does $\L_{X_\cdot^\mu}$ due to \eqref{u0}.

It remains to prove \eqref{ERT2}. To this end, we observe that \eqref{PDE}   and the It\^o formula yield
\beg{align*} \d \Theta^{\lambda,\mu}_t(X_t^\nu)&=\lambda \mathbf{u}^{\lambda,\mu}_t(X_t^\nu)\d t+(\{\nabla\Theta_t^{\lambda,\mu}\}\sigma_t)(X_t^\nu)\,\d W_t\\
&\qquad +[\{\nabla\mathbf{u}_t^{\lambda,\mu}\}(b^\nu_t-b^\mu_t)+ b_t^\nu-\bar b_t^\mu](X_t^\nu)\d t\\
&=[\lambda \mathbf{u}^{\lambda,\mu}_t+\{\nabla\Theta_t^{\lambda,\mu}\}(b^\nu_t-b^\mu_t)+\hat{b}^\mu_t](X_t^\nu)\d t+(\{\nabla\Theta_t^{\lambda,\mu}\}\sigma_t)(X_t^\nu)\,\d W_t.\end{align*}
 Combining this with  \eqref{E-X}  and applying the It\^o formula, we see that   $\eta_t:=\Theta^{\lambda,\mu}_t(X_t^\mu)-\Theta^{\lambda,\mu}_t(X_t^\nu)$ satisfies
\begin{equation*}\begin{split}
\d|\eta_t|^2
=&2\left<\eta_t,\lambda\mathbf{u}^{\lambda,\mu}_t(X_t^\mu)-\lambda\mathbf{u}^{\lambda,\mu}_t(X_t^\nu) +\hat{b}^\mu_t(X_t^\mu)-\hat{b}^\mu_t(X_t^\nu)\right\>\d t\\
&+2\left\<\eta_t,[(\{\nabla\Theta_t^{\lambda,\mu}\}\sigma_t)(X_t^\mu)-(\{\nabla\Theta_t^{\lambda,\mu}\}\sigma_t)(X_t^\nu)]\d W_t\right\>\\
&+\left\|(\{\nabla\Theta_t^{\lambda,\mu}\}\sigma_t)(X_t^\mu)-(\{\nabla\Theta_t^{\lambda,\mu}\}\sigma_t)(X_t^\nu)\right\|^2_{HS}\,\d t\\
&-2\left\<\eta_t, [\{\nabla\Theta_t^{\lambda,\mu}\}(b^\nu_t-b^\mu_t)](X_t^\nu)\right\>\d t.
\end{split}\end{equation*}
So, for any $m\ge 1$, it holds
\beq\label{NN1}\beg{split}
\d|\eta_t|^{2m}
 =\, &2m|\eta_t|^{2(m-1)}\left<\eta_t,\lambda\mathbf{u}^{\lambda,\mu}_t(X_t^\mu)-\lambda\mathbf{u}^{\lambda,\mu}_t(X_t^\nu)+ \hat{b}^\mu_t(X_t^\mu)-\hat{b}^\mu_t(X_t^\nu)\right\>\d t\\
&+2m|\eta_t|^{2(m-1)}\left\<\eta_t,[(\{\nabla\Theta_t^{\lambda,\mu}\}\sigma_t)(X_t^\mu)-(\{\nabla\Theta_t^{\lambda,\mu}\}\sigma_t)(X_t^\nu)]\d W_t\right\>\\
&+m|\eta_t|^{2(m-1)}\left\|(\{\nabla\Theta_t^{\lambda,\mu}\}\sigma_t)(X_t^\mu)-(\{\nabla\Theta_t^{\lambda,\mu}\}\sigma_t)(X_t^\nu)\right\|^2_{HS}\,\d t\\
&+2m(m-1) |\eta_t|^{2(m-2)}\left|[(\{\nabla\Theta_t^{\lambda,\mu}\}\sigma_t)(X_t^\mu)-(\{\nabla\Theta_t^{\lambda,\mu}\}\sigma_t)(X_t^\nu)]^\ast\eta_t \right|^2\d t\\
&-2m|\eta_t|^{2(m-1)}\left\<\eta_t, [\{\nabla\Theta_t^{\lambda,\mu}\}(b^\nu_t-b^\mu_t)](X_t^\nu)\right\>\d t.\end{split}\end{equation}
 By \eqref{u0} and \eqref{con},  we may find  a constant $c_0>0$ such that
\beq\label{XPP1}|\eta_t|^{2(m-1)} |\eta_t|\cdot|\lambda\mathbf{u}^{\lambda,\mu}_t(X_t^\mu)-\lambda\mathbf{u}^{\lambda,\mu}_t(X_t^\nu)+ \hat{b}^\mu_t(X_t^\mu)-\hat{b}^\mu_t(X_t^\nu)|\le c_0 |\eta_t|^{2m},\end{equation}
and
\beq\label{XPP3} \beg{split}
&|\eta_t|^{2(m-1)} |\eta_t| \cdot |[\{\nn\Theta_t^{\ll,\mu}\}(b_t^\nu-b_t^\mu)](X_t^\nu)|\\
&\le K_2 \|\nn\Theta^{\ll,\mu}\|_{\infty }|\eta_t|^{2(m-1)}|\eta_t|( \W_\theta(\mu_t,\nu_t)+\|\mu_t-\nu_t\|_{TV})\\
&\le c_0\big(|\eta_t|^{2m} +\W_\theta(\mu_t,\nu_t)^{2m}+\|\mu_t-\nu_t\|_{TV}^{2m}\big), \end{split}\end{equation}
According to \cite[(4.19)-(4.20)]{HW}, we arrive at
\beq\label{NNP}\d |\eta_t|^{2m} \le c_1 |\eta_t|^{2m} \d A_t + c_1 (\W_\theta(\mu_t,\nu_t)^{2m}+\|\mu_t-\nu_t\|_{TV}^{2m})\d t + \d M_t\end{equation}
for some constant $c_1>0$, a local martingale $M_t$,  and
$$
A_t:=\int_0^t\Big\{1 +\big(\scr M\big(\|\nn^2\Theta_s^{\ll,\mu}\|+\|\nn\si_s\|\big)(X^\mu_s)+  \scr M\big(\|\nn^2\Theta_s^{\ll,\mu}\|+\|\nn\si_s\|\big)(X_s^\nu)\big)^2\Big\}\d s.$$
Thanks to \cite[Theorem 3.1]{SQZ}, the Krylov estimate
  \beq\label{APP'}\begin{split}
 & \E\bigg[\int_s^t |f_r|( X^\mu_{r}) \d r\Big| \F_s\bigg]+  \E\bigg[\int_s^t |f_r|( X^\nu_{r}) \d r\Big| \F_s\bigg]\\
  &\le   C \left(\int_s^t \Big(\int_{\R^d} |f_r(x)|^p\d x\Big)^{\ff q p} \d r\right)^{\frac{1}{q}} ,\ 0\leq s<t\leq T.\end{split}\end{equation}
  holds.
As shown in \cite[Lemma 3.5]{XZ}, \eqref{APP'}, \eqref{HH2}, \eqref{u01} and $(A_\sigma)$ imply
$$\sup_{t\in [0,T]} \E\e^{ \delta A_t} = \E\e^{ \delta A_T}<\infty,\ \ \delta>0.$$
 By   \eqref{u0} and the stochastic Gronwall lemma (see  \cite[Lemma 3.8]{XZ}),    \eqref{NNP} with   $2m>\theta$ implies
 \begin{align*}&\{\W_\theta(\Phi_t^{\gg_1}(\mu),\Phi_t^{\gg_2}(\nu))\}^{2m}
    \le c_2 (\E |\eta_t|^\theta)^{\ff{2m}\theta} \\
 & \le c_3 (\E|X_0^\mu-X_0^\nu|^\theta)^{\ff{2m}\theta} + c_3\big(\E\e^{\ff{c_1\theta}{2m-\theta}A_T}\big)^{\ff{2m-\theta}{\theta} } \int_0^t(\W_\theta(\mu_s,\nu_s)^{2m}+\|\mu_s-\nu_s\|_{TV}^{2m})\d s
 \end{align*}
holds for all  $t\in [0,T]$ and some constants  $c_2,c_3>0.$   Therefore, \eqref{ENNB} holds for some constant $C>0$ and the proof is thus finished.

 \end{proof}

\section{Proof of Theorem \ref{T1.1}}

  Assume $(A_\si)$. According to \cite[Theorem 1.3]{Z2},  for any $\mu_\cdot\in \B([0,T];\scr P)$,   each of $(A_b)$ and $(A_b')$ implies the strong existence and uniqueness up to life time of the SDE \eqref{CSDE}. Moreover, it is standard that in both cases a solution of \eqref{CSDE} is non-explosive. So,    by Lemma \ref{YW}, the    strong well-posedness of  \eqref{E1} implies the weak well-posedness. Therefore, in the following we need only cosnider the strong solution.

To prove the strong well-posedness of \eqref{E1}, it suffices to find a constant $t_0\in (0,T]$ independent of $X_0$ such that in each of these two cases the SDE \eqref{E1} has  strong well-posedness up to time $t_0$. Indeed,
once this is confirmed, by considering the SDE from time $t_0$ we prove the same property up to time $(2t_0)\land T$. Repeating the procedure finite many times we derive the strong well-posedness.

Below we prove assertions (1) and (2) for strong solutions respectively.

(a) Let  $(A_b')$ hold.  Take $t_0= \min\{T,\ff 1 {K_1K_3^2}\} $  and consider the   space
$E_{t_0}:= \{\mu\in\B([0,t_0];\scr P):\mu_0=\gg\}$ equipped with the complete metric
$$\rr(\nu,\mu):= \sup_{t\in [0,t_0]} \|\nu_t-\mu_t\|_{TV}.$$  Then \eqref{ERT1} implies that $\Phi^\gg$ is a strictly contractive map on $E_{t_0}$, so that it has a unique fixed point, i.e. the equation
\beq\label{PHIE} \Phi_t^\gg(\mu)= \mu_t,\ \ t\in [0,t_0]\end{equation}  has a unique solution $\mu\in E_{t_0}.$ By \eqref{PHIE} and the definition of $\Phi^\gg$ we see that the unique solution of \eqref{CSDE} is a strong solution of \eqref{E1}.
On the other hand, $\mu_t:=\L_{X_t}$ for any strong solution to \eqref{E1} is a solution to \eqref{PHIE}, hence the uniqueness of \eqref{PHIE} implies that  of \eqref{E1}.

To prove \eqref{B1}, let $\mu_t= P_t^*\mu_0$ and $\nu_t =P_t^*\nu_0$. We have $P_t^*\mu_0= \Phi_t^{\mu_0}(\mu)$ and  $P_t^*\nu_0= \Phi_t^{\nu_0}(\nu).$  So,   \eqref{ERT1} with $\gg=\mu_0$ implies
\beq\label{B11} \|P_t^*\mu_0- \Phi_t^{\mu_0} (\nu)\|_{TV}^2\le \ff{K_1K_3^2} 4 \int_0^t \|P_s^*\mu_0-P_s^*\nu_0\|_{TV}^2\d s,\ \ t\in [0,T].\end{equation}
On the other hand, by the Markov property for the solution to \eqref{EW1} with $\nu$ replacing $\mu$, we have
$$\Phi_t^{\gg}(\nu)= \int_{\R^d} \Phi_t^{\dd_x} (\nu)\gg(\d x),\ \ \gg\in \scr P.$$ Combining this with $P_t^*\nu_0= \Phi_t^{\nu_0}(\nu)$, we obtain
\begin{align*}|\{\Phi_t^{\mu_0}(\nu)\}(A)- \{P_t^*\nu_0\}(A)|&=\bigg|\int_{\R^d}  \{\Phi_t^{\dd_x}(\nu)\}(A)(\mu_0-\nu_0)(\d x) \bigg|\\
 &\le \|\mu_0-\nu_0\|_{TV},\ \ A\in \B(\R^d).
\end{align*}
Hence,
\beq\label{B112} \|\Phi_t^{\mu_0}(\nu)-P_t^*\nu_0\|_{TV} \le \|\mu_0-\nu_0\|_{TV},\ \ t\in [0,T].\end{equation}
This together with \eqref{B11} yields
\beg{align*} &\|P_t^*\mu_0-P_t^*\nu_0\|_{TV}^2 \le 2 \|P_t^*\mu_0- \Phi_t^{\mu_0}(\nu)\|_{TV}^2 + 2 \| \Phi_t^{\mu_0}(\nu)-P_t^*\nu_0\|_{TV}^2\\
&\le 2 \|\mu_0-\nu_0\|_{TV}^2+ \ff{K_1K_3^2}2 \int_0^t\|P_s^*\mu_0-P_s^*\nu_0\|_{TV}^2\d s,\ \ t\in [0,T].\end{align*}
By Gronwall's lemma, this implies \eqref{B1}.

(b) Let $(A_b)$ hold and let $\gg=\L_{X_0}\in \scr P_\theta$. For any $\mu,\nu\in C([0,T],\scr P_\theta)$, \eqref{con} implies   \eqref{ESTM}.  By \eqref{ESTM}, \eqref{con}  and \eqref{ERT2} with $\gg_1=\gg_2=\gg$, we find  a constant $C>0$ such that
\beq\label{ERT3} \beg{split} & \{\|\Phi^\gg_t(\mu)-\Phi_t^\gg(\nu)\|_{TV}+  \W_\theta(\Phi_t^{\gg}(\mu), \Phi_t^{\gg}(\nu))\}^{2m} \\
&\le   C \int_0^t \big\{\|\mu_s-\nu_s\|_{TV}+\W_\theta(\mu_s,\nu_s)\big\}^{2m}\d s,\ \ t\in [0,T],\gg\in \scr P_\theta.\end{split}\end{equation}
Let $t_0= \ff 1 {2C}$. We consider the   space
$\tt E_{t_0}:= \{\mu\in C([0,t_0]; \scr P_\theta):\mu_0=\gg\}$ equipped with the complete metric
$$\tt \rr(\nu,\mu):= \sup_{t\in [0,t_0]} \{\|\nu_t-\mu_t\|_{TV}+\W_\theta(\nu_t,\mu_t)\}.$$
 Then $\Phi^\gg$ is strictly contractive in $\tt E_{t_0}$, so that the same argument in (a) proves the strong well-posedness of \eqref{E1} with $\L_{X_0}=\gg$ up to time $t_0$.

Let $\mu_t$ and $\nu_t$ be in (a). By \eqref{ERT3} with $\gg=\mu_0$ we obtain
\beq\label{B21} \begin{split} &\big\{\|P_t^*\mu_0-\Phi_t^{\mu_0}(\nu)\|_{TV}+\W_\theta(P_t^*\mu_0,\Phi_t^{\mu_0}(\nu))\big\}^{2m} \\
&\le C\int_0^t  \big\{\|P_s^*\mu_0-P_s^*\nu_0\|_{TV}+\W_\theta(P_s^*\mu_0,P_s^*\nu_0)\big\}^{2m}\d s,\ \ t\in [0,T].\end{split}\end{equation}
Next, taking   $\gg_1=\nu_0,\gg_2=\mu_0$ and   $\mu=\nu$ in \eqref{ERT2}, we derive
$$\big\{\W_\theta(P_t^*\nu_0, \Phi_t^{\mu_0}(\nu))\big\}^{2m}\le C \big\{\W_\theta(\mu_0,\nu_0)\big\}^{2m}.$$
Combining this with \eqref{B112} and \eqref{B21}, we find a constant $C'>0$ such that
\beg{align*}  &\big\{\|P_t^*\mu_0-P_t^*\nu_0\|_{TV}+\W_\theta(P_t^*\mu_0,P_t^*\nu_0)\big\}^{2m} \\
&\le  C'\big\{\|\mu_0-\nu_0\|_{TV}+\W_\theta(\mu_0,\nu_0)\big\}^{2m}  \\
&+C' \int_0^t  \big\{\|P_s^*\mu_0-P_s^*\nu_0\|_{TV}+\W_\theta(P_s^*\mu_0,P_s^*\nu_0)\big\}^{2m}\d s,\ \ t\in [0,T].\end{align*}
By Gronwall's lemma, this implies \eqref{B2} for some constant $c>0.$


\end{document}